

\input amstex
\documentstyle{amsppt}

\magnification=1200
\NoBlackBoxes

\TagsOnRight
\topmatter
\title
Some glueing formulas and the Van de Ven conjecture
\endtitle
\title Some glueing formulas for Spin polynomials and a proof of the Van de Ven
conjecture.
\footnotemark"*"
\endtitle
\footnotetext"*"{This is the translation of an article submitted to Izvestija
of
Russian Academy of Sciences, Ser.Math. While this translation was in progress
there appeared a paper [FQ2] with the proofs of the results announced in
[FQ1].}
\author Victor Pidstrigach
\endauthor
\address  Steklov Mathematical Institute, ul. Vavilova 42, Moscow, 117966,
GSP-1, Russia
\newline
e-mail: pidstrig\@alg.mian.su
\newline
 and
\newline
 SFB 343 Universit\"at Bielefeld, Postfach 100131, 33501 Bielefeld
\newline e-mail: pidstrig\@mathematik.uni-bielefeld.de
\endaddress

\abstract
We deduce some formulas for the spin invariants of the connected sum of an
arbitrary
4-manifold $X, b^+ _2(X) > 0 $ with $\overline {CP^2}$ in terms of the spin
invariants
of $X$ and apply this to prove the Van de Ven conjecture.
\endabstract
\endtopmatter


\define\sA{{\Cal A}} 
\define\sB{{\Cal B}} 
\define\sC{{\Cal C}} 
\define\sD{{\Cal D}} 
\define\sE{{\Cal E}} 
\define\sG{{\Cal G}} 
\define\sH{{\Cal H}} 
\define\sL{{\Cal L}} 
\define\sM{{\Cal M}} 
\define\sN{{\Cal N}} 
\define\sO{{\Cal O}} 
\define\sP{{\Cal P}} 
\define\sU{{\Cal U}}

\document

\head
\S 0. Introduction
\endhead

In this article we present some formulas for the spin-polynomial of the
connected sum of an arbitrary oriented smooth simply connected
four-manifold $X$ with $b_+ ^2(X) > 0$
with the orientation reversed complex projective plane. These formulas are
similar to ones for Donaldson's polynomials presented by Donaldson ([D])
and Morgan-Mrowka ([MM]) in the sense that they  do not involve evaluating
$ \nu$-classes for the polynomial of $X$. One can show that these formulas
also hold in the case of connected sums with arbitrary negative definite
manifolds. However the case of $\overline {CP^2}$ is easier from the technical
point of view.

As an application one has a proof of the Van de Ven conjecture [VdV]:
\proclaim{ Conjecture}
The Kodaira dimension of a simply connected algebraic surface is an
invariant of its smooth structure.
\endproclaim

This has been proved to be true for all cases except possibly nonminimal
surfaces
of general type (cf. [FM], [K], [OVdV], [P], [PT], [Q1], [Q2]). Therefore it
suffices
to prove that

\it
No surface of general type can be diffeomorphic to a rational one.
\rm

The proof is a reduction to the case of minimal surfaces of general type, which
already was
treated in [PT]. This reduction uses the glueing formulas
and results of [FM] describing diffeomorphisms
of rational surfaces .
In [FQ1]  a different approach to the Van de Ven conjecture is announced.
The authors use a glueing formula for Donaldson polynomials (and predict
our approach) and also explicitly formulate the result on diffeomorphism of
rational surfaces which we use.

\head
\S 1.Connected sum formulas.
\endhead

We shall start with some notations.

Let $X$  be a Riemannian, simply-connected 4-manifold with a metric $g$ and
$E \to X$ some
complex hermitian vector bundle of rank 2 with chern classes $c_1(E),
c_2(E)$. This vector bundle as well as the corresponding principal $U(2)$
-bundle
$P$ is uniquely determined by its chern classes. Denote also by $P_{ad}$ the
$PU(2)$-principal bundle $ P_{ad} = P/S^1$, and by $adE$ the corresponding
vector
bundle.

We shall consider the space $\sA_E$ of $U(2)$- connections on the bundle $E$
together with
the action of the gauge group $\sG_{U(2)}$ and the space of orbits of this
action $\sB_E$.
The space of $SO(3)$ -connections for the bundle $adE$, the gauge group, and
the orbit space, we denote by $\sA_{adE}, \sG_{SO(3)}, \sB_{adE} $
resp. We write $ \sA_{adE} ^*
\subset \sA_{adE}$ and $ \sB_{adE} ^* \subset \sB _{adE}$ for the set of
irreducible connections on $adE$ and the orbits of irreducible connections.

We shall also need the space of orbits of framed connections
$ \widetilde{\sB_{adE}}$, that is the orbit
space of the space of connections by the subgroup
$ \sG^0 _{G} $ of the gauge group acting trivially on the fibre over some
fixed point $p$.

Let
$$ \sA_{\lambda} = \{ A \in \sA_E {~}|{~}\frac12 trF_A = \lambda \} $$
for some closed differential 2-form $\lambda $ on $X$.
There is an isomorphism:
$$ \sA_{\lambda}/\sG_{U(2)} = \sA_{adE} / \sG_{(SO(3)=PU(2))} .\tag{1.1}$$

Let $C$ be some $Spin^{\Bbb C}$-structure on the manifold $X$.
This means that one fixes an integral lift $C$ of the second Stiefel-Whitney
class
$w_2 (X)$ of $X$, and this lift defines a pair of complex hermitian rank 2
vector bundles
$W^+, W^- $ subject to the following conditions:
$$
T^*_{\Bbb C} X = Hom_{\Bbb C} (W^+, W^-),
{~~~} c_1 ( \Lambda ^2 W^{\pm} ) = C.
$$
A choice of a connection $ \nabla$ on the line bundle $\Lambda ^2 W^+ $
gives rise to a Dirac operator on $X$:

$$ \Cal D^{C,\nabla} : \Gamma (W^+ ) \rightarrow \Gamma (W^- ) .$$

Coupling it with connections $a \in \Cal A_E$ then gives a family of Fredholm
operators over $\sA_E$:

$$ \Cal D_a ^{C,\nabla} : \Gamma (W^+ \otimes E ) \rightarrow \Gamma
(W^- \otimes  E) ,$$
which is equivariant with respect to the gauge group action.
Using the isomorphism (1.1) we consider it as a family over $\sB _{adE}$.

\it We shall always suppose that \rm
$$ ind( \Cal D_a ^{C,\nabla} ) = dim ker
(\Cal D_a ^{C,\nabla}) - dim coker (\Cal D_a ^{C,\nabla}) \le 0 .
$$

Let $ \Bbb S ( \Gamma (W^+ \otimes E)) $ denote the unit sphere in the Hilbert
space $ \Gamma (W^+ \otimes E)  $, and
$ \Bbb P ( \Gamma (W^+ \otimes E)) $  the projectivisation of this space.
There is a tautological bundle
$
\sO_{   \Bbb P ( \Gamma (W^+ \otimes E))} (1)
$
over this projective space,
and the coupled Dirac operator can be interpreted as a section of the bundle
$$
\Gamma (W^- \otimes E) \otimes \sO_{\Bbb P ( \Gamma (W^+ \otimes E))} (1).
$$

The space
$$ \sP = \sA_{\lambda} \times_{\sG_{U(2)}} \Bbb S ( \Gamma (W^+ \otimes E) )
\tag{1.2}$$
is the base of the following Hilbert bundle
$$
\sH_C = \sA_{\lambda} \times_{\sG_{U(2)}} Tot,
$$
where $Tot$ stands for the total space of the bundle
$ \Gamma (W^- \otimes E) \otimes \sO_{\Bbb P ( \Gamma (W^+ \otimes E)} (1)$.
One can see that the family of coupled Dirac operators gives a section $s_D$
of the bundle
$\sH_C \to \sP$
(details in  1.1.26,  [PT]).
In the same way one gets a section $s_{asd}$ of the Hilbert bundle
$\sH_{+}$ associated to the principal bundle $\sA_{adE} ^* \to
\sB_{adE} ^*$ and a natural representation of the gauge group $\sG_{adE}$
in the Hilbert space $ \Omega ^+ _2 (adE) $.

\definition{Definition}
The moduli space of pairs
$ \Cal {MP}^{g,\nabla}(X, p_1 (E), C+c_1(E)) $
is the space of zeroes of the section $s_D \oplus \pi ^* s_{asd}$ where
$\pi: \sP^* \to \sB^*$ is the natural projection.

\enddefinition

Using (1.1) one can give a point of $ \Cal {MP}^{g,\nabla}(X, p_1(E), C+c_1(E))
$
as a pair $(a, \left< s \right>), a \in \sM ^g _{asd} (X, p_1(E), w_2(E))$ and
$\left< s \right>$ is a point of the projective space
$ \Bbb P(\Gamma (W^+ \otimes E))$
given by a complex line spanned by a nonzero vector spinor $s$.
There is an exact sequence which describes the tangent bundle of
moduli of pairs at a point $(a, \left< s \right>)$ :
$$ 0 \to T_{a,s} \Cal {MP}^{g,\nabla}
\to (ker \Cal D_a)/\left< s \right> \oplus H^1_a (adE) \overset i \to
\to coker \Cal D_a
\tag{1.3}$$
$$ i(s', \alpha) = \varpi(\alpha * s) ,$$
(here $*$ is a combination of spinor multiplication and multiplication
of a matrix by a vector, $\varpi$ is the orthogonal projection
to the space of harmonic vector spinors $coker \Cal D_a$).
Therefore  $ T_{a,s} \Cal {MP}^{g,\nabla}$ is the kernel of the linearisation
of our asd- and Dirac- equations and $coker (i) $ is its cokernel.
One can globalise this exact sequence to a fibre of $\pi$ as
an exact sequence of bundles over
$\Bbb P (ker \Cal D_a ^{C,\nabla} )$ :
$$
0 \to T \Cal {MP}^{g,\nabla}_{|\pi ^{-1}a}
\to T \Bbb P (ker \Cal D_a ^{C,\nabla} ) \oplus H^1_a (adE)
\overset i \to \to  coker \Cal D_a
\otimes \Cal O_{ \Bbb P (ker \Cal D_a ^{C,\nabla} )} (1) \to 0 .
$$
An orientation of the moduli space of pairs is defined by the orientation
of the moduli space of asd connection via the exact sequence (1.3).

A model for a neighborhood of the point
$(a,\left< s \right>)$ in the space
$\Cal {MP}^{g,\nabla} (X   )$
is given by the zero set of some smooth section of the bundle

$$
\CD
(coker (i)) \oplus H_{a}^2(adE) \otimes \sO_{\Bbb P (ker(\Cal
D_{a}^{C,\nabla}))} \\
@VVV \\
\Bbb P (ker(\Cal D_{a}^{C,\nabla})) \times H^1_{a} (adE)
\endCD
$$
(cf. [PT], ch.1).

This moduli space is a source of various invariants of the smooth structure
of $X$ as described in [PT] ( numerical invariants) and [T] (polynomial
invariants). It is a smooth manifold provided $b_2 ^+ (X) >0$ and the metric
$g$ and the connection $\nabla$ both are generic.
For $b_2^+(X)$ odd we use polynomial invariants
$$ \gamma^X _{p_1, C+c_1(E)} (\Sigma_i) = \#  \{ \cap V_{ \Sigma_i} \cap
\Cal {PM} (X, p_1(E), C+c_1(E))\}, $$
where $ \Sigma _i \in  H_2(X) $,
and the  $V_{ \Sigma_i}$ are codimension 2 submanifolds of $\sP $,
defined as the zero sets of certain sections $s_{\Sigma}$ of certain
complex line bundles $ \sL _{ \Sigma_i}$  over $\sP$ (cf. [T], [D]).
Such a section is defined once we have chosen a smooth Riemann surface
representing
the class $\Sigma$.
We shall use the same notation for these representatives: $\Sigma_i$.
If $ b^+ _2 >1 $ then these polynomials do not depend on the choice
of the metric $g$ and the connection  $ \nabla$ on the line bundle $\Lambda ^2
W^+ $ .
If $b^+ = 1$ and $p_1 \ge -7 $  or $p_1 = -8 , w_2 = (c_1) \pmod 2 \ne 0 $
then there is a chamber structure on the space of periods $\sP(X)$ of metrics
and the polynomials depends only on the chamber, not on the metric itself (cf.
[BP]). We
shall use the notation  $\gamma^X _{p_1, C+c_1(E), \sC} (\Sigma_i)$ to
show the dependence on the chamber $\sC$.
In either case  as a function of the cohomology classes $c_2(E), c_1(E), C$ the
polynomials
depend in fact only on the  combinations
$C+c_1(E), p_1(adE) = c_1(E) ^2 - 4c_2(E)$ of these classes.

The connected sum $X \# \overline {CP^2}$ of a manifold $X$ and a projective
plane with reversed
orientation $\overline {CP^2}$ can be viewed as a union
$$
(X \setminus D_1) \sqcup (S^3 \times [-t/2,t/2]) \sqcup (\overline {CP^2}
\setminus D_2)
$$
where $D_1 \subset X, D_2 \subset \overline {CP^2}$ are small disks around
some points $p \in X, q \in \overline {CP^2}$, $\partial D_i = S^3$. Define a
metric
$g_t$ on $X \# \overline {CP^2}$  which is generic when restricted to $X -D_1$
and
$\overline {CP^2} - D_2$ and
a product of $S^3$ with radius $r$ and a standard metric on the interval $[
-t/2, t/2]$
for  some interior of the cylinder.

When $t$ tends to infinity one gets two manifolds - $X$ and $\overline {CP^2}$
joined by an "infinite"cylinder or, in a conformally equivalent metric,
two open manifolds
$(X \setminus p) $ and $(\overline {CP^2}\setminus q)$.

Let $e$ is a generator of $H^2(\overline {CP^2})$.
Fix some $U(2)$ vector bundle $E$ and $Spin^{\Bbb C}$ structure on the
connected sum $X \# \overline {CP^2}$,
such that $(C+c_1(E))_{ \vert \overline {CP^2}}  = ke$ with $ -3<k<3$.

It is known that if one takes a sequence of metrics $g_{t_i}, t_i
\to \infty$ and a sequence of
connections $a_i$ on $E$ such that $a_i$ is antiselfdual (asd) with respect to
$g_{t_i}$ then there exist asd-connections $a_X, a_{\overline {CP^2}} $ on
$X,\overline {CP^2} $ resp., and a finite collection
of points $x_1, ..., x_m$  such that after gauge transformation some
subsequence of $a_i$ converges weakly over compact
subsets of $X \cup \overline {CP^2} - \{ p,q,x_1, ..., x_m\}$ to
$a_X,a_{\overline {CP^2}} $ and the curvature
densities $\vert F_{a_i} \vert$ converge to
$F_{a_X} +F_{a_{\overline {CP^2}}} + \sum_i
8 \pi ^2 \delta_{x_i}$ (cf. ch.7 of [DK]).

The following lemma describes the case when there is a nontrivial element
of the kernel of the coupled Dirac operator $D_{a_i}$ (positive harmonic
vector-spinor) for any asd-connection $a_i$.

\proclaim{ Lemma 1.1}  Take $E, C, a_i$ as above.
If for any asd-connection $a_i$  with $i>>0$ there exists a nonzero
positive harmonic vector-spinor $s_i$ (which we shall suppose to be of unit
norm) , then $a_X$ also has nonzero positive harmonic vector-spinor $s_X$.
Moreower a subsequence of $s_i$ converges to $s_X$ on any compact subset of
$X-p$ in $C ^{\infty}$ topology.

\endproclaim

\demo{Proof}
We shall use the norm $L^{8/3}$ on the space $\Gamma(W^+
\otimes E)$ of positive vector spinors which has following property:
if $g' = e^{2f} g$ is a conformal change of a metric then
$$
s \in ker (\Cal D_a ^g) \Leftrightarrow s' = e^{- \frac32 f} s \in ker (\Cal
D_a ^{g'})
$$
and
$$
\Vert s \Vert _{L^{8/3} (X,g)} = \Vert s' \Vert  _{L^{8/3} (X,g')}.
$$
This gives a possibility to derive estimates in an arbitrary conformal model of
the
connected sum, identifying vector spinors as above.

Let $S_t$ be a scalar cuvature of $ X \# \overline {CP^2}, g_t$. It is easy to
see that it is
uniformly bounded from above for all $t$. Let $ \nabla _{a_i} $ be a connection
on the
bundle $W^+ \otimes E$ built from $a_i$ and $\nabla$.
It follows from the Weitzenb\"ock formula
$$
(\Cal D_{a_i})^* \Cal D_{a_i} s_i = \nabla _{a_i} ^* \nabla _{a_i} s_i -
F_{a_i} ^+ s_i + 1/4 S_t \cdot s_i
+F_{\nabla} s_i ,
$$
that if $\sigma$ the maximum of the norm of the sum of curvatures $ 1/4 S_t
+F_{\nabla}$
then there is the pointwise inequality:
$$
( \nabla_{a_i} ^* \nabla _{a_i} s_i, s_i)  \le
(\sigma) ( s_i, s_i) .
$$
Being combined with the standart estimate of the l.h.s. through Laplasian
$\Delta$ of the function
$\vert s_i \vert$
$$
( \nabla_{a_i} ^* \nabla _{a_i} s_i, s_i) \ge \vert s_i \vert \cdot \Delta
(\vert s_i \vert),
$$
which follows from the Kato inequality, it gives then
$$
\Delta (\vert s_i \vert) \le \sigma \vert s_i \vert {~~~~}{~~~~} \text{and}
{~~~}{~~~~~}
\Vert \Delta (\vert s_i \vert) \Vert _{L ^{8/3}} \le \sigma \Vert s_i \Vert _{L
^{8/3}}.
$$
Since $\Delta$ is an elliptic operator one has an a priori inequality for
functions on $ X \# \overline {CP^2}$:
$$ \Vert s_i \Vert _{L^{8/3} _2}  \le const \Vert \Delta s_i \Vert
_{L^{8/3}} + const  \Vert s_i \Vert_{L^{8/3}} $$
with $const$ independent of $t$ and $a_i$.
The last two estimates give
$$
\Vert s_i \Vert _{L^{8/3} _2}  \le const(1+ \sigma) \Vert s_i \Vert _{L^{8/3}}.
$$
Applying the Sobolev embedding theorem for functions one gets
$$\Vert s_i \Vert _{L^{8}} \le const'  \Vert s_i \Vert _{L_2 ^{8/3}} \le
const'' \Vert s_i \Vert _{L^{8/3}}
$$
with $const''$ independent of $t,a_i$. With this estimate in mind we now look
how
far the truncated spinor is from being harmonic.

There are some fixed $\theta_1, \theta_2$ with  $-1/2< \theta_1, \theta_2  <1/2
$
such that for a subsequence of $a_i$ restricted to
the subcylinder $S^3 \times [ \theta_1 t, \theta_2 t] $ the curvature density
$\vert F_{a_i} \vert $ is small. Thus we can apply the theorem of K. Uhlenbeck
on the existence of local Coulomb gauge with small connection matrix over the
subcylinder.

Now let $ \psi_i $ be some cut-off function with support on
$X \cup S^3 \times [-t/ 2, \theta_2 t]$,
$supp (d \psi_i) \subset S^3 \times [ \theta_1 t, \theta_2 t] $  and
$\Vert d \psi_i \Vert _{L^{4}} \le \epsilon (t_i)$, where
$\epsilon (t_i) \to 0 $ when $t_i  \to \infty $.
We assume for simplicity of notations that there is no bubbling points on $X$ -
the modification
to this case consists of cutting off also in small neighborhood of the bubbling
points and has
been treated in [DKr, Lemma 7.1.24]. Therefore we can assume that there is an
estimate for
the supremum of $(a_i - a_X)$:
$$
sup (\vert a_i -a_X \vert _{\vert supp \psi _i}) < \epsilon ' (t), {~}{~}{~}
\epsilon ' (t) \to 0 {~~~} \text{as} {~~~} t \to \infty .
$$

This provides the estimates
$$
\Vert \Cal D_{a_X} \psi_i s_i \Vert _{L^{8/3}(X)} \le
$$
$$
\le \Vert \Cal D_{a_i} s_i \Vert _{L^{8/3}(X \# \overline {CP^2})} +
\Vert d \psi_i \cdot s_i \Vert _{L^{8/3}(X \# \overline {CP^2})} +
sup \vert a_i - a_X \vert _{\vert supp \psi_i}
\Vert  s_i \Vert _{L^{8/3}(X \# \overline {CP^2})} \le
$$
$$
\le \Vert d \psi_i  \Vert _{L^{4}(X)} \Vert s_i \Vert _{L^{8}(X\# \overline
{CP^2})} +
sup \vert a_i - a_X \vert _{\vert supp \psi_i} \Vert  s_i \Vert _{L^{8/3}(X\#
\overline {CP^2})} \le
$$
$$
\le (const'' \epsilon (t) + \epsilon' (t) ) \Vert  s_i \Vert _{L^{8/3}(X\#
\overline {CP^2})}
= \epsilon_1(t_i) \Vert  s_i \Vert _{L^{8/3}(X\# \overline {CP^2})} =
\epsilon_1 (t_i), \tag{1.4}
$$
where $\epsilon_1(t_i)$ tends to $0$ if $t_i$ tends to infinity.

On the other hand applying these estimates to the cut-off function
$(1- \psi_i)$ provides a similar inequality:
$$
\Vert \Cal D_{a_{\overline {CP^2}}} (1-\psi_i) s_i \Vert _{L^{8/3}} \le
\epsilon_2.
$$
But the auxilary Lemma 1.2 shows the first eigenvalue of the operator
 $\Cal D_{a_{\overline {CP^2}}}$
is positive and uniformly estimated from below by some positive constant $c$.
Therefore there is an inequality
$$
\Vert (1-\psi) s_i \Vert _{L^{8/3}} \le \epsilon_2/c.
$$
So the norm of $(1-\psi_i) s_i$ is arbitrarily small and hence
the norm of $\psi_i s_i$ is arbitrarily close to $1$.
Now the first statement of the Lemma follows from inequality (1.4) and the
variational
characterisation of the first eigenvalue.
The second statement follows from standart bootstrapping argument using the
equation
$$
\Cal D_{(a_X)_{\vert K}}  (s_i)_{\vert K} = (a_i - a_X)_{\vert K} * (s_i)
_{\vert K}
$$
restricted to the compact subset $K$ of $X-p$.
\enddemo

\proclaim{Remark}
If one has
$$
a_i \in V_{\Sigma}, \Sigma \in
H_2(X) \subset H_2(X \# \overline {CP^2})
$$
and no bubbling occures at $\Sigma$, i.e. none of the $x_i$ lies on  $ \Sigma$
then for the limit connection it follows
$$
a_X \in V_{\Sigma},
$$
since this depends only on the restriction of the connection to
a neighborhood of the Riemaniann surface $ \Sigma$ where there is $C^{\infty}$
convergence of connections of our sequence.
\endproclaim

\proclaim{Lemma 1.2} Let $F$ be a $U(2)$ vector bundle on
$\overline {CP^2}$ and $C$ a $Spin^{\Bbb C}$-structure on $\overline {CP^2}$
such that
$C+c_1(E) = ke, -3<k<3$, where $e$ is a generator of
$H^2(\overline {CP^2})$.
There exists a neighborhood $U$ of the Fubini-Study metric on $\overline
{CP^2}$ and
some connection in the product of the space of metrics on $\overline {CP^2}$
and the space of connections on the line bundle $L=\Lambda ^2 W^{\pm},
c_1(L) = ke$, such that
for any asd-connection on $adF$ and any point  $(g, \nabla)  \in U$, one has

$$ ker \Cal D_a ^{g,\nabla}  = 0.$$
Moreover, the first eigenvalue  $\mu_a$ of the Laplasian
$ ( \Cal D_a ^{g,\nabla} )^* \Cal D_a ^{g,\nabla} $ satisfies
$$ \mu_a > c $$
for some positive constant  $c$ independent of the point in $U$ and $a$.

\endproclaim
\demo{Proof}
The proof is based on the Weitzenb\"ock formula:
$$ (\Cal D_a)^* \Cal D_a s = (\nabla _a) ^* \nabla _a s - F_a ^+ s + 1/4 S
\cdot s +F_{\nabla} s , \tag{1.5}$$
where  $S$ is the scalar curvature of $\overline {CP^2}$, $ \nabla _a $ the
connection on the bundle $W^+ \otimes F$ built from  $a$ and $\nabla$.
One has $F_a ^+ = 0$ , and example 2.1.2 [H] gives that for the Fubini-Study
metric, some connection $\nabla_0$ on
$\Lambda ^2 W^{\pm}$ and some positive number $c$  there is an
estimate of $1/4S +  F_{\nabla_0}$ as an endomorphism of vector-spinors:
$$
1/4S +  F_{\nabla_0} > 2 c.
$$
Let $U$ be given by the condition
$(g, \nabla)  \in U$ if $1/4S +  F_{\nabla} >  c$.
Taking the scalar product of (1.5) with $s$ and integrating over $\overline
{CP^2}$ gives
the required estimate.

\enddemo

Therefore to each point of the space
$ \Cal {MP}^{g_i,\nabla} (X \# \overline {CP^2}, C+c_1 , p_1) $ for $i \to
\infty$
one can assign a point of the product
$$
{\Cal {MP}}^{g,\nabla} (X, (C+c_1)_{\vert X}, k) \times
{\Cal {M}^g}( \overline {CP^2}, (w_2)_{ \vert \overline {CP^2}}, m),
$$
$w_2 = c_1 \pmod 2$, at least for one pair of negative integers
$k, m; k+m \ge p_1$. Conversely, there is a procedure for glueing
two asd-connections $a_X, a_{\overline {CP^2}},$ defined on the summands,
which gives an asd connection on the connected sum
and one shows that any asd-connection
on a connected sum with long enough tube is obtained in this way (cf. [DKr]).
In the following we extend this construction to the case of the moduli space of
pairs.

To start with we need glueing data for connections, i.e. isometries
$$
\tau: T_p X \to T_q \overline {CP^2},
\rho: (adE_X)_p \to (adE_{\overline {CP^2}})_q,
$$
which for given asd-connections $a_X, a_{\overline {CP^2}}$ on $ X$and $
\overline {CP^2}$
resp. give an $SO(3)$-asd-connection $a = I(a_X, a_{\overline {CP^2}},
\rho) $  on $X \#\overline {CP^2}$ with appropriate choices of some
other parameters . Varying $a_X, a_{\overline {CP^2}}$ in some
small neighborhood, and varying the glueing parameter $\rho$ in the space of
all glueing
parameters $Gl$ we get a smooth map to
some open set of asd connections on the connected sum.

In order to glue vector-spinors
one has to glue the corresponding vector bundles of vector-spinors with a
unitary map
$ \sigma$:
$$
\sigma: (W_X ^+ \otimes E_X)_p \to (W_{\overline {CP^2}} \otimes E_{\overline
{CP^2}})_q .
$$
The isometries $\tau $ and $\rho$ give $\sigma$  up to a unitary scaling, which
may be interpreted as a glueing parameter for glueing  connections on the
determinant
line bundles
$$
(det(W_X ^+ \otimes E_X))_p \to (det(W_{\overline {CP^2}} \otimes E_{\overline
{CP^2}}))_q.
$$
This scaling is in fact inessential since different glueing parameters for
$U(1)$- connections give gauge equivalent connections provided $X$ is simply
connected.
Taking a cut-off function $\psi$ as in Lemma 1.1 one considers $ \psi s$ as a
section
of the glued bundle $ (W^+ \otimes E)_{X \# \overline {CP^2}}$ and it turns out
that
the point $\left( a =  I(a_X, a_{\overline {CP^2}}, \rho), \left< \psi s
\right> \right)$
is close to  a point
of  $\Cal {MP}^{g,\nabla}(X \# \overline {CP^2}, C+c_1, p_1)$. In fact we can
deform it
to a point of the moduli space. The next lemma describes a neighborhood of this
deformation
in the space $ \Cal {MP}^{g_i,\nabla} ( X \# \overline {CP^2}, C+c_1, p_1)$,
as the zero set of some map defined on some neighborhood in the space
$$
\widetilde {\Cal {MP}}^{g,\nabla} (X, (C+c_1)_{\vert X},-4 k) \times _{PU(2)}
\widetilde {\Cal {M}^g} ( \overline {CP^2}, (w_2)_{ \vert \overline {CP^2}}  ,
p_1 +4k)
$$
of the preimage of the pair
$$
((a_X,\left< s \right>) ; a_{\overline {CP^2}}) \in {\Cal {MP}}^{g,\nabla} (X,
(C+c_1)_{\vert X},-4 k) \times
{\Cal {M}^g}( \overline {CP^2}, (w_2)_{ \vert \overline {CP^2}} , p_1 +4k)
$$
under the natural projection
$$
\widetilde {\Cal {MP}}^{g,\nabla} (X, (C+c_1)_{\vert X}, -4k) \times _{PU(2)}
\widetilde {\Cal {M}^g} ( \overline {CP^2}, (w_2)_{ \vert \overline {CP^2}} ,
p_1 +4 k)
\to
$$
$$
{\Cal {MP}}^{g,\nabla} (X, (C+c_1)_{\vert X}, -4k) \times
{\Cal {M}^g}( \overline {CP^2}, (w_2)_{ \vert \overline {CP^2}} , p_1 +4 k).
$$

\proclaim{Lemma 1.3} For the point
$(a_X, \left< s \right>)$ in the space ${\Cal {MP}}^{g,\nabla} (X,
(C+c_1)_{\vert X}, -4k)$
and $a_{\overline {CP^2}} \in  {\Cal {M}^g}( \overline {CP^2}, (w_2)_{ \vert
\overline {CP^2}} , p_1 +4 k)$
there is a neighborhood in the space
$ \Cal {MP}^{g_i,\nabla} (X \# \overline {CP^2}, C+c_1, p_1) $
which is diffeomorphic to the zero set of a smooth map $ \Phi$ defined on some
neighborhood
$$
\sU \subset
T_{a,s} \Cal {MP}^{g,C,\nabla}
\times H_{a_{\overline {CP^2}}}^1(adE) \times Gl
$$
of $ \{ 0 \} \times \{ 0 \} \times Gl$:

$$\CD
\sU @>\Phi>>
coker (i)  \oplus (H_{a_X}^2(adE)
\oplus H_{a_{\overline {CP^2}}}^2(adE))
\oplus coker \Cal D_{a_{\overline {CP^2}}}
\endCD
$$
where $coker(i)$ is defined in (1.3).
\endproclaim

\demo{Proof}
We assume for simplicity that $H^2 _{a_X} = H^2  _{\overline {CP^2}} = 0$.
Fix some glueing parameter $\rho'$ and a connection $a' = I(a_X  , a_{\overline
{CP^2}}, \rho')$.
Consider a neighborhood in $\Cal {M}^g (X \# \overline {CP^2}, (c_1) \pmod2,
p_1)$
of the subspace  $ \cup _{\rho} I(a_X  , a_{\overline {CP^2}},  \rho)$
consisting
of all connections represented in the form
$I(a_X +\delta a , a_{\overline {CP^2}},  \rho)$ subject to an estimate of
$\delta I(a_X +\delta a , a_{\overline {CP^2}},  \rho)
= I(a_X +\delta a , a_{\overline {CP^2}},  \rho) - I(a_X, a_{\overline {CP^2}},
 \rho')$:
$$
\Vert \delta I(a_X +\delta a , a_{\overline {CP^2}},  \rho) \Vert _{L^4} \le
\epsilon
$$
with $\epsilon$ to be fixed in the proof.
Fix lifts
$$
\xi _X : coker (\sD_X) \to \Gamma (W^- \otimes E_X), {~~~}
\xi _{\overline {CP^2}} : coker \sD _{a_{\overline {CP^2}}} \to
\Gamma (W^- \otimes E_{\overline {CP^2}})
$$
of the obvious projections $\pi_X, \pi_{\overline {CP^2}}$ onto the cokernels
with all section in the images $ im \xi _X, im \xi_{\overline {CP^2}}$
supported away
from the points $p \in X, q \in \overline {CP^2}$. There is a linear map
$$
\xi _{X \# \overline {CP^2}} =  \xi _X \oplus  \xi_{\overline {CP^2}}
:coker \sD _{a_X} \oplus coker \sD _{a_{\overline {CP^2}}}
\to \Gamma (W^- \otimes E_{\bar {X \# \overline {CP^2}}}),
$$
which is a lift of the projection
$ \pi (\cdot) = \pi_X ( \psi \cdot) + \pi_{\overline {CP^2}} ((1- \psi)
\cdot)$,
where $\psi$ is the cut-off function as in Lemma 1.1.
Let $s$ be a non-zero element of the kernel of $\sD_{a_X} $. Then we can make
$ \Vert \sD_{a'} (\psi s) \Vert \le \epsilon (t)$ with $\epsilon (t)  \to 0$
as $ t \to \infty$ for suitable choices as above.
Consider the map $\phi $, which assignes to a triple consisting of an
asd-connection $a_X + \delta a_X$ on
$X$, a vector spinor $\delta s$ on the connected sum $X \# \overline {CP^2}$
and an element $ h \in coker \sD _{a_X} \oplus coker \sD _{a_{\overline
{CP^2}}} $
the vector spinor on the connected sum:
$$
\phi (a_X + \delta a_X, h, \delta s) =
\sD _{I(a_X +\delta a_X , a_{\overline {CP^2}},  \rho)} ( \psi s + \delta s)
+\xi_{X \# \overline {CP^2}}(h).
$$
The tangent map of $\phi$ at the point $ (a_X, 0, 0)$
with respect to $\delta s, h$ is given by
$$
(\sD_{a'} \oplus \xi) (h, \delta s) = \sD _{a'} (\delta s) + \xi_{X \#
\overline {CP^2}} (h),
$$
and the application of the standart technique shows that it has a right inverse
$$
S \oplus \pi : \sD_{a'} S (\eta)  +
\xi_{X \# \overline {CP^2}} \pi (\eta ) = \eta
$$
with the estimate of the norm of $S \oplus \pi$ independent of $t$ (compare
7.2.14, 7.2.18 of  [DKr]).
Since by definition $pr _2 = \xi_{X \# \overline {CP^2}} \pi$ is a projector,
$pr _1 = \sD_{a'} S$ is also a projector and one can also assume that
$ Im S$ is transversal to $ ker \sD_{a'} $. This means that $S (\eta) = S (
\eta_1)$,
where $ \eta_i = pr_i (\eta)$.
Although $pr_i$ are not orthogonal projectors one can estimate their norm by,
say,
$\Vert pr_i \Vert \le 2$.
Now we shall look for the solution of the equation
$ \phi (a_X + \delta a_X, h, \delta s) = 0$
in the form $(h, \delta s) = (\pi (\eta), S(\eta)) = (\pi (\eta_2), S(
\eta_1))$:
$$
\sD _{I(a_X + \delta a_X, a_{\overline {CP^2}}, \rho} (\psi s - S (\eta)) + \xi
(\pi (\eta)) =
$$
$$
=\sD_{a'} (\psi s) + \eta_1 + \eta_2 +  \delta I(a_X +\delta a_X , a_{\overline
{CP^2}},  \rho) * S(\eta)
+ \delta I(a_X +\delta a_X , a_{\overline {CP^2}},  \rho) *( \psi s) = 0,
\tag (1.6).
$$
Following  system of two equations is equivalent to (1.6):
$$
pr _1 (\sD_{a'} (\psi s)) + \eta_1 +
pr_1( \delta I(a_X +\delta a_X , a_{\overline {CP^2}},  \rho)* (\psi s)) +
pr_1( \delta I(a_X +\delta a_X , a_{\overline {CP^2}},  \rho) * S(\eta_1)) = 0,
\tag{1.7}$$
$$
pr_2(\sD_{a'} (\psi s)) + \eta_2 +
pr_2( \delta I(a_X +\delta a_X , a_{\overline {CP^2}},  \rho)* (\psi s)) +
pr_2(\delta I(a_X +\delta a_X , a_{\overline {CP^2}},  \rho) * S(\eta_1)) = 0.
\tag{1.8}$$
Equation (1.7) can be rewritten in the form:
$$
(1+A)(\eta_1) = - pr _1 (\sD_{a'} (\psi s)) -
pr_1( \delta I(a_X +\delta a_X , a_{\overline {CP^2}},  \rho)* (\psi s)
\tag{1.9}
$$
with the norm of $ A( \eta_1) = pr_1( \delta I(a_X +\delta a_X , a_{\overline
{CP^2}},  \rho) * S(\eta_1))$
estimated by an inequality:
$$
\Vert A( \eta_1) \Vert _{L^{8/3}} \le
2 \Vert \delta I(a_X +\delta a_X , a_{\overline {CP^2}}, \rho) \Vert _{L^4}
\Vert S( \eta_1)  \Vert_{L^8} \le
2 \epsilon \Vert S \Vert \Vert \eta_1 \Vert_{L^{8/3}} \le \Vert \eta_1
\Vert_{L^{8/3}} /2
$$
provided $ \epsilon \le \Vert S \Vert /4$.
Therefore the operator $1+A$ is invertible and $ \Vert (1+A) ^{-1} \Vert \le
2$.
This gives the existence and uniqueness of the solution $\eta_1$ of (1.9)
together with an estimate:
$$
\Vert \eta_1 \Vert_{L^{8/3}} =
\Vert (1+A)^{-1} (- pr _1 (\sD_{a'} (\psi s)) -
pr_1( \delta I(a_X +\delta a_X , a_{\overline {CP^2}},  \rho)* (\psi s))
\Vert_{L^{8/3}} \le
$$
$$
\le 2 (\Vert pr_1(\sD_{a'} (\psi s)) \Vert_{L^{8/3}} +
\Vert pr_1 (\delta I(a_X +\delta a_X , a_{\overline {CP^2}},  \rho)* (\psi s))
\Vert_{L^{8/3}}) \le
$$
$$
\le 4 (\Vert \sD_{a'} (\psi s) \Vert_{L^{8/3}} +
\Vert  \delta I(a_X +\delta a_X , a_{\overline {CP^2}},  \rho)* (\psi s)
\Vert_{L^{8/3}}) \le
$$
$$
\le 4( \epsilon (t) + \Vert \psi \delta I(a_X +\delta a_X , a_{\overline
{CP^2}},  \rho) \Vert_{L^4}
\Vert s \Vert_{L^8} ) \le 4( \epsilon (t) + const_1 \epsilon )
\tag {1.10}$$
(we use here the estimate $\Vert s \Vert _{L^{8}} \le const' \Vert s \Vert
_{L_2 ^{8/3}} \le
const" \Vert s \Vert _{L^{8/3}} $ from Lemma 1).
Let $\sN$ be the image of some splitting $j$ of the exact sequence (1.3).
Now consider the map $\tau : \sN \to coker \sD _{a_X} \oplus coker \sD
_{a_{\overline {CP^2}}}$
given by
$$
\delta a_X \mapsto
\pi_X ( \delta I(a_X +\delta a_X , a_{\overline {CP^2}},  \rho)* (\psi s)) +
\pi_X (\delta I(a_X +\delta a_X , a_{\overline {CP^2}},  \rho) * S(\eta_1)).
$$
In order to get the proper obstruction space
we shall prove that the image of the composition $pr_{im( i)} \tau $
containes an open ball centered at the origin, where $pr_{im(i)}$ is
the projection to the image of the map $i$ defined in (1.3).
Indeed, this follows from estimates:
$$
\Vert \pi_X ( \delta I(a_X +\delta a_X , a_{\overline {CP^2}},  \rho)* (\psi
s)) +
\pi_X (\delta I(a_X +\delta a_X , a_{\overline {CP^2}},  \rho) * S(\eta_1))
\Vert_{L^{8/3}} \ge
$$
$$
\ge \Vert \pi_X ((\delta a_X) *(\psi s) \Vert_{L^{8/3}} /2  -
\Vert \delta I(a_X +\delta a_X , a_{\overline {CP^2}}, \rho) \Vert_{L^4} \Vert
S \Vert
\cdot \Vert \eta_1 \Vert_{L^{8/3}}
\ge
$$
$$
\ge \Vert \pi_X ( (\delta a_X) * s) \Vert_{L^{8/3}} /4 -
2 \Vert \delta a_X \Vert_{L^4} \Vert S \Vert \cdot 4 (\epsilon(t) +const_1
\epsilon)=
$$
$$
= \frac 14 \Vert  i_{\vert \sN } ( \delta a_X )) \Vert_{L^{8/3}} -
8 \Vert S \Vert (\epsilon(t) +const_1 \epsilon) \Vert (\delta a_X) \Vert_{L^4}
\ge
$$
$$
\ge (const_2 - 8 \Vert S \Vert (\epsilon(t) +const_1 \epsilon)) \Vert ( \delta
a_X) \Vert_{L^4},
\tag{1.11}
$$
where the first inequality uses the estimate
$$
\Vert \delta I(a_X +\delta a_X , a_{\overline {CP^2}},  \rho) - \psi \cdot
\delta a_X
\Vert _{L^4} \le const \cdot  e^{-2t}
\tag{1.12}
$$
(cf. [DKr, (7.2.37)]). The second inequality in (1.11) uses (1.12), (1.10) and
the estimate of the norm of $ \psi s$ from below given in
Lemma 1.1. Since $i_{\vert \sN}$ is a linearisation of $\pi_X ( (\delta a_X) *
s)$
(cf.1.3) and
$i_{\vert \sN}$ is a monomorphism by the definition of $\sN$, there
is a constant $const_2$ such that for small enough $\epsilon$
$$
\Vert \delta a_X \Vert \le \epsilon \Rightarrow \Vert \pi_X ( (\delta a_X) * s)
\Vert \ge
\frac 12
\Vert  i_{\vert \sN } ( \delta a_X )) \Vert_{L^{8/3}} \ge const_2 \Vert (
\delta a_X) \Vert_{L^4}
$$
(here we identify $\delta a_X$ with the corresponding element in $H^1 _{a_X}$).
This gives the last inequality in (1.11), and  the statement about the
composition
$pr_{im (i)} \tau $ follows, provided $\epsilon(t)$ and $ \epsilon$ are small
enough.
Therefore the subspace $ im( \tau)  + coker (i) \subset coker (\sD_{a_X}) $
contains some metric ball $B$, and
there exist $\epsilon_0$ such that
$$
\epsilon (t) \le \epsilon_0 \Rightarrow
\pi (  \sD_{a'} (\psi s)) \in B \oplus coker \sD _{a_{\overline {CP^2}}}
\subset coker \sD _{a_X}
\oplus coker \sD _{a_{\overline {CP^2}}}.
$$
Therefore the obstruction space can be identified with $coker (i) \oplus coker
\sD _{a_{\overline {CP^2}}} $.
The Lemma follows.

\enddemo

Denote by $\epsilon _{\Omega} $ the constant in
the generalisation of Uhlenbeck's local Coulomb gauge theorem to any
strongly simply connected domain $\Omega$ (Prop.4.4.10 of [DKr]).

\proclaim{Corollary 1.4}
There exists a real positive number $t_0$ such that for $t>t_0$ there exists
 a covering $\{ \sU_i \}$ of the moduli space
$ \Cal {MP}^{g, \nabla} (X \# \overline {CP^2}, C+c_1, p_1) $
such that for subsets
$$
\widetilde{\Cal {MP}^{\varrho}} ( X, (C+c_1) _{\vert X}, i)
\subset \widetilde{\Cal {MP}} ( X, (C+c_1) _{\vert X}, i))
$$
and
$$
\widetilde{\Cal {M}^{\varrho}} ( \overline {CP^2}, (w_2) _{\vert \overline
{CP^2}} , p_1+4i )
\subset \widetilde{\Cal {M}} ( \overline {CP^2}, (w_2) _{\vert \overline
{CP^2}} , p_1+4i )
$$
given by the condition
$$
\Vert F_{a_{\vert D^4}} \Vert _{L_2 (D^4)}  \le \epsilon _{D_{\varrho} ^4}
$$
where $D^4 _{\varrho}$ is a small disk around $p \in X$ , $q \in \overline
{CP^2}$,
resp., of radius $\varrho$, one has

$$
\sU_i \subset \widetilde{\Cal {MP}}^{r exp(-t)} ( X, (C+c_1) _{\vert X}, -4i)
\times
_{PU(2)}
 \widetilde{\Cal {M}}^{r exp(-t)} ( \overline {CP^2}, (w_2) _{\vert \overline
{CP^2}}
 , p_1+4i )
$$
$$
\widetilde{\Cal {MP}^{r }} ( X, (C+c_1) _{\vert X}, -4i) \times _{PU(2)}
 \widetilde{\Cal {M}^{r }} ( \overline {CP^2}, (w_2) _{\vert \overline {CP^2}}
, p_1+4i ) \subset \sU_i.
$$
\endproclaim

\demo{Proof}
Take a subdivision of the tube $S^3 \times [-t/2, t/2]$ into $N$ pieces
$$
\Omega _l =  S^3 \times [-t/2+ lt/N,- t/2 +(l+1)t/N],  l = 0,...,N-1.
$$
We shall use a cut off function $\psi$ with support of $d \psi$ in one of
$\Omega _l$
with an estimate $ \Vert d \psi \Vert _{L^4} \le (10 (t/N)^3 ) ^{-1/4}$.
Take $ N-2  > max \{ -p_1 / \epsilon _ F ,10 (\epsilon_0)^{-4} t^{-3} \}$,
where
$\epsilon _F \le \epsilon _ {\Omega_l}$ and $ \epsilon_0 $ is the constant of
Lemma 1.3.
This provides the existence for any asd-connection $a$ with first Pontrjagin
class $p_1$ of at least one $0 <l< N-1$ for which
$$
\Vert F_{a_{\vert \Omega_l}} \Vert _{L_2 }
\le \epsilon _ F, {~~~~~}
\Vert \sD _a ( \psi s) \Vert_{L^{8/3}(\Omega _l)} \le
\Vert d \psi \Vert _{L^4 (\Omega _l)} \Vert s_{\vert \Omega_l} \Vert \le
(10 (t/N)^3)^{-1/4}  \Vert s_{\vert \Omega_l} \Vert \le
$$
$$
(10 (t/N)^3)^{-1/4} N^{-1} = (10 t^3 N)^{-1/4} \le \epsilon_0
\tag {1.13}$$
Now let $ \sU _i $ be the open set of those connections $a$ in
$ \Cal {MP}^{g, \nabla} (X \# \overline {CP^2}, C+c_1, p_1) $
for which for at least one integer $l$ (1.13) is satisfied, and the $L_2$-norm
of the curvature of $a$ restricted to $X_l =X \cup S^3 \times [-t/2,-t/2+
lt/N]$
is close to $8 \pi ^2 i$.
Then by the Uhlenbeck theorem there is a trivialisation of
$E_{\vert \Omega _l}$ in which the
connection matrix is small in $L_{k} ^2$ and therefore one can show that norms
of
$
F^+ _{ \psi a}, F^+ _{(1- \psi) a} , \sD_{ a} (\psi s)
$
are small. By Theorem 7.2.41 of [DKr]
$ a = I (a_X, a_{\overline {CP^2}}, \rho) $,  $\Vert \psi a - a_X \Vert _{L^q}
\le const \epsilon _F $.
If $l<N-1$, $X_l$ is a subset of fixed compact subspace in $X-p$ and by
Lemma 1.1 one has for some $ s_X \in ker \sD_{a_X}$
$$
\Vert  (s - s_X)_{\vert X_l} \Vert _{L^{8/3} _1 (X_l)} \to 0
$$
as $t \to \infty$. Symmetrically, if $l>0$ then
$$
\Vert s _{\vert (X \# \overline {CP^2} - X_l)} \Vert _{L^{8/3} _1 (X \#
\overline {CP^2} -X_l)} \to 0.
$$
This means that
$$ \Vert  s - \psi  s_X \Vert _{L^{8/3} _1 (X \#\overline {CP^2})}
\le
\Vert  (s - s_X)_{\vert X_l} \Vert _{L^{8/3} _1 (X_l)} +
$$
$$
+\Vert s _{\vert (X \# \overline {CP^2} - X_l)} \Vert _{L^{8/3} _1 (X \#
\overline {CP^2} -X_l)}
+ \Vert s \Vert _ {L^{8/3} _1 (\Omega_l)} + \Vert \psi \Vert _{L^2 _2 (\Omega
_l)}
\Vert s_X \Vert _ {L^{8/3} _1 (\Omega_l)} \to 0 $$
as $t \to \infty$, so it follows that varying slightly $s_X$ in $ker \sD_{a_X}$
(what does not change $S_a$) we get
$s - \psi  s_X = S_a ( \eta) $ with $ \eta = \sD_a ( s - \psi  s_X) =
(\sD_{a_X} + (\psi a - \psi a_X) )(s_X)$
small enough to be in the set of solutions of  (1.7 -1.8).
\enddemo

Now we shall consider some particular cases in which Lemma 1.3
will provide sufficient information about the moduli space of
the connected sum $X \# \overline {CP^2}$ in order to compute the
spin-polynomials of the connected sum in terms of those of $X$.
Consider a spin polynomial of degree $d$ on $X$:
$$
\gamma^X _{k, \omega} (\Sigma_i) = \#  \{ \cap V_{ \Sigma_i } \cap
\Cal {PM} (X, k, \omega) \},
$$
which is defined only if
$$
2d = -(3/2)k  -3(1+b_2^+(X))+(\omega)^2/2 -sign(X)/2 -2
$$
for some negative integer $k$ and some  2-dimensional cohomology class $\omega$
on $X$.

Let $e$ be a generator of the second cohomology group $H_2( \overline {CP^2})$:
$ e \in H_2( \overline {CP^2}) \subset H_2(X \# \overline {CP^2})$.
Take a bundle $E$ and a $Spin ^{\Bbb C }$ -structure $C$ on the
connected sum
$X \# \overline {CP^2}$
given by one of the following three choices.

1. $p_1(E) = k,  C+c_1(E)=\omega+(\pm e)$;

2. $ p_1(E) = k +1 , C+c_1(E)=\omega+(\pm 2e)$;

3. $ p_1(E) = k +1 , C+c_1(E)=\omega $.

By a dimension count it follows that in the first and second case one has
$$
\dim \left( \cap V_{ \Sigma_i} \cap \Cal {PM} (X \# \overline {CP^2}, p_1(E),
C+c_1(E)) \right) =
0 ,
$$
and in the third -
$$
\dim \left( \cap V_{ \Sigma_i} \cap \Cal {PM} ( X \# \overline {CP^2}, p_1(E),
C+c_1(E)) \right) =
2.
$$

Therefore in the first and second case the  polynomial
$$
\gamma^{X \# \overline {CP^2}} _{p_1, C+c_1(E)} (\Sigma_i) =
\# \{ \cap V_{ \Sigma_i} \cap
\Cal {PM} (X \# \overline {CP^2}   , p_1(E), C+c_1(E)) \} ,
$$
is defined and in the third the polynomial
$$
\gamma^{X \# \overline {CP^2}} _{p_1, C+c_1(E)} (\Sigma_i, e) =
\# \{ \cap V_{ \Sigma_i}  \cap V_e \cap
\Cal {PM} (X \# \overline {CP^2}   , p_1(E), C+c_1(E))\}
$$
is defined.

Let $a_j$ be a sequence of connections in $\cap V_{ \Sigma_i}$
which are asd in the metric $g_j$ resp.
In each of the three cases by the remark to Lemma 1.1 the limit connection
$a_X$ on $X$ belongs to the moduli space
$$\cap V_{ \Sigma_i} \cap \Cal {PM} (X, k+4m, \omega), $$
for some $m \ge 0 $,
of the virtual dimension $-8m$. Since we took the metric on $X$ and the
connection on the bundle  $\Lambda ^2 W^{\pm}$ to be generic,
and since $b_2^+(X) > 0$ there are no reducible asd-connections and the
moduli spaces of negative virtual dimension are empty
(cf. [PT], ch.1, sec.3).
Therefore the only remaining case is  $m = 0$.

{}From the inequality
$$
p_1(a_X) + p_1(a_{\overline {CP^2}}) \ge p_1 (E)
$$
it follows that $p_1(a_{\overline {CP^2}}) = 0 $ for the first case and
$p_1(a_{\overline {CP^2}}) = -1 $ for the second and the third.
Moduli space of asd-connections on $\overline {CP^2}$ for a generic metric
with these Pontrjagin numbers
has negative virtual dimension and therefore contain only reducible
connections. It follows that $a_{\overline {CP^2}}$ is the trivial connection
in
the first case, and the only reducible one with $p_1 = -1$ in the second.

In this situation Lemma 1.3 provides a satisfactory description of
the (compact) manifold
$$
Z = \cap \tilde{V} _{ \Sigma_i} \cap \widetilde{\Cal {PM}} (X, k, \omega)
\times _{PU(2)} \tilde M (\overline {CP^2}, w_2( a_{\overline {CP^2}}), p_1(
a_{\overline {CP^2}}))
$$
and the bundle $\Xi$ on $Z$ which is a tensor product of the lift of the
$PU(2)$ -equivariant bundle
$$
U(2) \times _{S^1 \times S^1} coker \sD_{  a_{\overline {CP^2}}}
\to \tilde M (\overline {CP^2}, w_2( a_{\overline {CP^2}}), p_1( a_{\overline
{CP^2}}))
= U(2)/(S^1 \times S^1)
$$
via obvious projection
$$ \CD
 \cap \tilde{V} _{ \Sigma_i} \cap \widetilde{\Cal {PM}} (X, k, \omega)
\times _{PU(2)} \tilde M (\overline {CP^2}, w_2( a_{\overline {CP^2}}), p_1(
a_{\overline {CP^2}})) \\
@VVV \\
\tilde M (\overline {CP^2}, w_2( a_{\overline {CP^2}}), p_1( a_{\overline
{CP^2}}))
\endCD
$$
and tautological line bundle on $\cap \tilde{V} _{ \Sigma_i} \cap
\widetilde{\Cal {PM}} (X, k, \omega)$
such that our, cut down moduli space
$$
\cap V_{ \Sigma_i} \cap \Cal {PM} (X \# \overline {CP^2}, p_1(E), C+c_1(E))
$$
is embedded in $Z$ as the zero set of some generic section of $\Xi$.
This gives a possibility of computing the polynomial for $X \# \overline
{CP^2}$
in topological  terms.

In the first case our manifold $Z$ is parametrised as a family of glued
connections of the type $a_X \# \theta$ with
$(a_X, <s_X>) \in \cap V_{ \Sigma_i} \cap \Cal {PM} (X , p_1(E), C+c_1(E)) $.
Since a connected sum with the
trivial connection $\theta$ does not need any glueing parameters, one has
$ Z= \cap V_{ \Sigma_i} \cap \Cal {PM} (X, k, \omega)$. Index of the coupled
Dirac operator for $\theta$ is zero and hence $ \Xi = 0$.

In the second case the reducible connection on $\overline {CP^2}$ has
the stabilizer group $S^1$ and therefore
$$
Z= \cap \tilde{V} _{ \Sigma_i} \cap \widetilde{\Cal {PM}}
(X, k, \omega)/ S^1 = \cup \Bbb P^1.
$$
For each component $\Bbb P^1$ the restriction to it of the bundle $\Xi$ is
$$
\left( (U(2) \times U(2)) \times _{S^1 \times (S^1 \times S^1)}
(\Bbb C \otimes coker \sD_{  a_{\overline {CP^2}}}) \right) /U(2)
$$
where $S^1 \times (S^1 \times S^1)$ is the subgroup given as the product
of the center subgroup $S^1 \subset U(2)$ and the maximal torus
$(S^1 \times S^1) \subset U(2)$. The representation of this subgroup is
the product of the generating linear one of $S^1$ and the natural
representation
of the stabiliser subgroup $(S^1 \times S^1)$ of the reducible connection
$a_{\overline {CP^2}} = \lambda_1 \oplus \lambda_2$ which is the sum
of representations $\oplus_i coker \sD _{\lambda_i}$, each $S^1$ acting
on its summand with degree $\pm 1$.
For the
specified $Spin ^{\Bbb C}$ -structure on $\overline {CP^2}$ this
cokernel is one dimensional and therefore
$$
\Xi_{\vert \Bbb P ^1}  = U(2) \times _{S^1 \times S^1} coker \sD_{
a_{\overline {CP^2}}} =
\sO _{\Bbb P^1} (1).
$$

In the third case the space $Z$ is the same, but the index of the reducible
connection on $\overline {CP^2}$ is zero, the Dirac operator coupled
to this connection is an isomorphism and the
bundle arises as a line bundle $ \sL_e, c_1( \sL_e) =  \mu(e) $ :
$$
\Xi = \sL _{ e} = \sO _{\Bbb P^1} (-1).
$$

Therefore the following theorem is proved:

\proclaim{Theorem }
Let $k \in \Bbb Z, k \le 0;  \omega \in H^2(X \Bbb Z), \Sigma_1, ..., \Sigma _d
\in H_2(X, \Bbb Z)$
such that
$$ 2d = -3/2 k -3(1+b_2^+(X))+(\omega)^2/2 -sign(X) -2 \ge 0.$$
Then
$$ \gamma^{X \# \overline {CP^2}} _{k, \omega \pm e} (\Sigma_1, ..., \Sigma_d)
=
\pm \gamma^X _{k, \omega} (\Sigma_1, ..., \Sigma_d)  ,$$

$$ \gamma^{X \# \overline {CP^2}} _{k+1, \omega \pm 2e} (\Sigma_1, ...,
\Sigma_d) =
\pm \gamma^X _{k, \omega} (\Sigma_1, ..., \Sigma_d) ,$$

$$ \gamma^{X \# \overline {CP^2}} _{k+1, \omega} (\Sigma_i, ..., \Sigma_d, e) =
\pm \gamma^X _{k, \omega} (\Sigma_1, ..., \Sigma_d). $$

\endproclaim

\proclaim{Remarks}

1.Signes in the r.h.s. of the formulas are determined by the choice of
orientation.

2.Considering the case $b^+ _2 (X) = 1 $ one has to take into account the
chamber
structure i.e. the fact that stretching the neck of the connected sum
$X \# \overline {CP^2}$ provides periods of our metric tending to the
hyperplane
of orthogonals to $e$ in the limit. However in the case of Spin-polynomials
this
hyperplane is not a wall by the Proposition 1.4.2 of [PT] and the negative
definiteness
of the intersection form on $H^2( \overline {CP^2}).$
Therefore any chamber in  the period space $\sP (X)$  for metrics on $X$ is
embedded in a unique chamber in the period space $ \sP (X \# \overline
{CP^2})$.
\endproclaim

\head
\S 2.Application
\endhead

Let $S_r$ and $S_g$ be a rational surface and a surface of general type with
$p_g(S_g) = 0$ resp. We assume the rational surface to be a good generic
surface
(cf. [FM, 1.2.1]).
Denote by $p : S_g \to S_{mg} $ the map to the minimal model of $S_g$ and by
$\{ L_i \}$ a collection  of
exceptional fibres of this map.
For $e \in H^2(S_r, \Bbb Z),  (e)^2 = -1$ the wall $W^e$ will be a hyperplane
of orthogonals
to  $e$ in $H^2(S_r, \Bbb Z)$. One says that the wall is ordinary if the
restriction of the intersection
form to the wall is odd, and extraordinary otherwise.
For $ e, (e)^2 = -1$, the reflection $ R_e $ in the wall $W^e$ is defined by
$ R_e (x) = x+ 2(x.e)e$. It follows from [FM,  Prop.3.2.4] that if $e$ is
represented by an embedded (-1)-sphere,
i.e. there is an embedded sphere $S^2$
with the fundamental class Poincare dual to $e, (e)^2 = -1$, there is an
orientation preserving self-diffeomorphism
of $S_r$ inducing
$R_e$. This is the case in particular for any $e \in \sE(S_r)$ where $\sE(S_r)$
is a
set of all cohomology classes represented by exceptional (-1)-curves on $S_r$.
We shall use the notions of a P-cell and a super P-cell as defined in [FM,
ch2].

\proclaim{Lemma 2.1}
\footnote{Proved also in [FQ2] as Theorem 1.7.
The arguments of our proof are a minor modification of the arguments in [FM]
and are of
course  similar to those in [FQ2].}
Suppose $f$ is a hypothetical diffeomorphism
$$
S_r \to S_g.
$$
There exists a self-diffeomorphism $\psi$ of $S_r$ and a collection of
exceptional
curves $E_i$ on $S_r$ such that  for all (-1)-curves $L_i$ on $S_g$ we have
$$
\psi ^* f^* l_i = e_i
$$
with $l_i = [L_i] , e_i = [E_i] $ denoting the 2-cohomology classes in
$H_2(S_g) , H_2(S_r)$ represented
by $L_i, E_i$ resp.
\endproclaim
\demo{Proof}
Let $\Cal K$  be the K\"ahler cone of $S_r$. Take any polarisation $H \in \Cal
K$.
Consider  the projection $\prod_i \frac{(1+ R_{f^* l_i})}2 $  onto the subspace
$\cap_i W^{f^*l_i}$
and consider the point
$$
p= \prod_i \frac{(1+ R_{f^* l_i})}2 H
$$
in this subspace. Since the cohomology class $f^* l_i$ is represented by a
(-1)-sphere
$f^{-1}( L_i)$ the reflection  $ R_{f^* l_i}$ is induced
by a diffeomorphism. Thus does any reflection of the type
$ R_I = \prod_{i \in I}  R_{f^* l_i} $.
Therefore by Thm.10 of [FM] it follows that $R_I$ leave the super P-cell $\Bbb
S ( \Cal K)$ invariant
and $R_I (H) \in \Bbb S ( \Cal K) $ . By
the convexity of the supercell (cf. Prop. 2.5.6 of [FM]) one has therefore
$p \in \Bbb S ( \Cal K)$.

Take any P-cell $\Cal P$ in $\Bbb S ( \Cal K)$ which  contains the point $ p
\in \Cal P$ under
consideration. By the definition of the super P-cell there exists an isometry
$\psi$ in  the subgroup of isometries generated by reflections in the set $\Cal
F _{ \Cal K}$
of ordinary walls of $\Cal K$ such that
$ \psi ( \Cal P ) = \Cal K $.
But these latter reflections are induced by diffeomorphisms since all such
walls are given by (-1)-curves on $S_r$ (cf. Prop.3.6
and Thm10B of [FM]). Therefore, $\psi $ is induced by a diffeomorphism.

Now, any wall $W^e$ for $ e = \psi ^* f^* l_i$
contains the point $ \psi^* (p) $ of the P-cell $ \Cal K$. In fact $W^e $
is a wall of the P-cell $ \Cal K$. To establish this fact one needs to show
that the
intersection $W^e \cap \Cal K$ containes an open subset of $W^e $. By
Proposition 3.4 of [FM]
$$
\Cal K = \{ y \in \Bbb H_+ | y \cdot (-K_{S_r}) \ge 0 ,{~~~} y \cdot h \ge 0
\text{ for all classes} {~~~}
h \in \sE(S_r) \}.
$$
Therefore
$$
W^e \cap \Cal K = \{ y \in \Bbb H_+ | y \cdot pr(-K_{S_r}) \ge 0 , {~~~} y
\cdot pr(h)\ge 0
\text{ for all classes} {~~~}  h \in \sE (S_r) \},
$$
where $pr$ is the projection to $W^e$. Let $E_1, ..., E_k$ be the set of all
exceptional
curves satisfying $ q \in W^{[E_i]}$ and assume that $\pm e \ne E_i$ (otherwise
we are done).
It follows from Lemma 1.11 of ch.2 of [FM] that
$ [E_i] = pr [E_i] \in W^e $ is a set of pairwise orthogonal elements. For some
neighborhood
$U_q$ of $q \in W^e$ one has
$$
U_q \cap W^e \cap \Cal K = \{ y \in W^e \cap U_q | y \cdot pr(-K_{S_r}) \ge 0 ,
{~~~}y \cdot pr([E_i] )\ge 0 \},
$$
and this intersection is an open set in $W^e$  if the elements $ pr (-K_{S_r}),
[E_i] \in W^e$ are
linearly independent. But if this is not the case one has $ K_{S_r} = \pm a e$.
This is impossible
since $W^{K_{S_r}}$ is even and $e$ defines an ordinary wall unless $S_g$ is
homeomorphic to
$(CP^1 \times CP^1) \# \overline{CP^2}$ and in this case $K_{S_g} ^2 = 7$ is
not a complete square.

So $W^{\psi ^* f^* l_i} $ is a collection of ordinary walls of the P-cell $\Cal
K$ with a nonempty
intersection $\cap W^{\psi ^* f^* l_i} = \psi^* f^* H^2(S_{gm})$,
and for any such collection one has $\psi ^* f^* l_i = [E_i] $ , as follows
from
Prop3.6, ch.2 of [FM].
\enddemo

\proclaim{Corollary 2.2}
No surface of general type can be diffeomorphic to a rational one.
\endproclaim

\demo{Proof}
Assume that  there exists a diffeomorphism $f : S_r \to S_g$. By Lemma 2.1 one
can also assume that there exists a  blow-down $S_r \to \tilde{S}$ such that
$f^* (H^2(S_{mg}) = H^2( \tilde {S})$.
It follows from [PT] that there are classes
$p_1 \in H^4(S_{mg}), \beta \in H^2(S_{mg}), p_1 > -8$  and a chamber $\sC
\subset \sP (S_{mg})$
such that
$$
\gamma _{p_1, \beta, \sC} ^{S_{mg}} \ne 0, {~~~}
\gamma _{p_1, f^* \beta, f^* \sC} ^{\tilde S} = 0.
$$

Now by Theorem  we have a contradiction
$$
0 \ne \gamma _{p_1, \beta, \sC} ^{S_{mg}} =
\gamma _{p_1, \beta + \sum l_i, \sC} ^{S_{g}} =
\pm \gamma _{p_1, f^* \beta + \sum e_i, f^* \sC} ^{ S_r} =
\pm \gamma _{p_1, f^* \beta, f^* \sC} ^{\tilde S} = 0.
$$
\enddemo

{\bf Acknowledgment.} The author is grateful to Frank-Olaf Schreyer for the
kind invitation and
the Mathematical Institute of the University of
Bayreuth for their hospitality during research on this project. The author is
also grateful to
Stefan Bauer for his kind invitation and
to the SFB 343 at Bielefeld University for the support during the translation
of
this text. Thanks to Rogier Brussee and Manfred Lehn for help with the
translation of this article
and once more to Rogier Brussee for pointing out the wrong integer factor in
one of the glueing formulas.

\enddocument


\Refs
\widestnumber\key{OVdV}

\ref\key BP
\by S.Bauer and V.Pidstrigatch
\paper Spin polynomial invariants for Dolgachev survaces
\jour alg-geom/9311008
\endref

\ref\key D
\by  S.Donaldson
\paper Polynomial invariants for smooth four-manifolds
\jour Topology
\vol 29
\yr 1990
\pages 257-315
\endref

\ref\key DKr
\by S.Donaldson and P.Kronheimer
\book The geometry of four-manifolds
\bookinfo Oxford Mathematical Monographs
\publ Oxford University Press
\yr 1990
\endref

\ref\key FM
\by R.Friedman and J.W.Morgan
\paper On the diffeomerphism types of certain algebraic surfaces I
\jour J. Diff. Geom.
\vol 27
\yr 1988
\pages 297-369
\endref

\ref\key FQ1
\by R.Friedman and Z.Qin
\paper The smooth invariance of the Kodaira dimension of a complex surface
\jour Mathematical Research Letters 1
\yr 1994
\pages 369-376
\endref

\ref\key FQ2
\by R.Friedman and Z.Qin
\paper On complex surfaces diffeomorphic to rational surfaces
\jour alg-geom/9404010
\endref

\ref\key H
\by N.J.Hitchin
\paper Harmonic spinors
\jour Advances in Math.
\vol 14
\yr 1974
\pages 1-54
\endref

\ref\key K
\by D.Kotschick
\paper On manifolds homeomorphic to $ CP^2 \# 8 \overline {CP^2}$
\jour Invent. Math.
\vol 95
\yr 1989
\pages 591-600
\endref

\ref\key MM
\by J.Morgan and T.Mrowka
\paper A note on Donaldson's polynomial invariants
\jour Int. Math. Research Notices
\vol 10
\yr 1992
\pages 222-230
\endref

\ref\key OVdV
\by C.Okonek and A. Van de Ven
\paper $ \Gamma$-invariants associated to $ \Bbb PU(2)$-bundles and the
differentiable
structure of Barlow surface.
\jour Invent.Math.
\vol 95
\yr 1989
\pages 602-614
\endref

\ref\key P
\by V.Pidstrigatch
\paper Deformation of instanton surfaces
\jour Math USSR Izv.
\vol 38
\yr 1992
\pages 313-331
\endref

\ref\key PT
\by V.Pidstrigatch and A.Tyurin
\paper Invariants of the smooth structure of an algebraic surface arising from
the Dirac operator
\jour Russian Acad. Sci. Izv. Math.
\vol 40
\yr 1993
\pages 267-351
\endref

\ref\key T
\by A.Tyurin
\paper The Spin-polynomial invariants of the smooth structures of algebraic
surface
\jour Izv. Russ. Acad. Nauk, Ser. Mat.
\vol 57
\yr 1993
\lang Russian
\transl\nofrills English transl. in
\jour Mathematica Gottingensis
\vol 6
\yr 1993
\transl\nofrills and to appear in
\jour Russian Acad. Sci. Izv. Math.
\yr 1994
\endref


\ref\key Q1
\by Z.Qin
\paper Complex structures on certain differentiable 4-manifolds
\jour Topology
\vol 32
\yr 1993
\pages 551-566
\endref

\ref\key Q2
\by Z.Qin
\paper On smooth structures of potential surfaces of general type
homeomorphic to rational surfaces
\jour Invent.Math.
\vol 113
\yr 1993
\pages 163-175
\endref

\ref\key VdV
\by A. Van de Ven
\paper On the differentiable structure of certain algebraic surfaces
\jour Sem. Bourbaki
\vol 667
\yr Juin 1986
\pages
\endref



\endRefs
\end {document}